\documentclass{article}

\usepackage{amsmath}
\usepackage{amsthm}
\usepackage{alltt}
\usepackage{hyperref}
\usepackage{algorithm}
\usepackage{algpseudocode}
\usepackage{graphicx}

\newtheorem{theorem}{Theorem}[section]

\numberwithin{equation}{section}

\newtheorem{lemma}{Lemma}
\numberwithin{lemma}{section}

\newcommand{\multsum}[1]
{ \mathop{\sum\cdots\sum}_{
   \scriptstyle k_1=0\cdots\  k_{#1}=0
   }^{\infty\ \ \cdots\ \ \infty}
}

\title{Computation of $P(n,m)$, the Number of\\ Integer Partitions of $n$ into Exactly $m$ Parts}
\author{M.J. Kronenburg}
\date{}

\begin{document}

\maketitle

\begin{abstract}
Two algorithms for computing $P(n,m)$, the number of integer partitions
of $n$ into exactly $m$ parts, are described, and using a combination of these two algorithms,
the resulting algorithm is $O(n^{3/2})$. The second algorithm uses a list
of $P(n)$, the number of integer partitions of $n$, which is cached
and therefore needs to be computed only once. Computing this list is also $O(n^{3/2})$.
With these algorithms also $Q(n,m)$, the number of integer partitions of $n$ into exactly
$m$ distinct parts, and a list of $Q(n)$, the number of integer partitions of $n$
into distinct parts, can be computed in $O(n^{3/2})$.
A list of $P(n,1)..P(n,n)$ and $P(m,m)..P(n,m)$ can be computed
in $O(n^2)$. A computer algebra program is listed implementing these algorithms,
and some timings of this program are provided.

\end{abstract}


\section{Definitions and Basic Identities}

Let the coefficient of a power series be defined as:
\begin{equation}
 [q^n] \sum_{k=0}^{\infty} a_k q^k = a_n
\end{equation}
Let $P(n)$ be the number of integer partitions of $n$, and let $P(n,m)$ be the
number of integer partitions of $n$ into exactly $m$ parts.
\begin{theorem}\label{pnq}
\begin{equation}
 P(n) = [q^n] \frac{1}{\prod_{j=1}^n (1-q^j)}
\end{equation}
\end{theorem}
\begin{proof}
Using the geometric series:
\begin{equation}
 [q^n] \frac{1}{\prod_{j=1}^n (1-q^j)} = [q^n] \prod_{j=1}^n \sum_{k_j=0}^{\infty} q^{jk_j}
  = [q^n] \multsum{n} q^{\sum_{j=1}^n jk_j} = P(n) 
\end{equation}
\end{proof}
\begin{theorem}\label{pnmq}
\begin{equation}
 P(n,m) = [q^n] \frac{q^m}{\prod_{j=1}^m (1-q^j)}
\end{equation}
\end{theorem}
\begin{proof}
\begin{equation}
 [q^n] \frac{q^m}{\prod_{j=1}^m (1-q^j)} = [q^n] q^m \prod_{j=1}^m \sum_{k_j=0}^{\infty} q^{jk_j} 
 = [q^n] \multsum{m} q^{m+\sum_{j=1}^m jk_j} = P(n,m)
\end{equation}
where $P(n,m)$ is the number of integer partitions of $n$ with greatest part equal to $m$.
By conjugation of Ferrer diagrams this is also the number of integer partitions of $n$
into exactly $m$ parts \cite{AE04}.
\end{proof}
\begin{theorem}\label{recur}
\cite{WPPP}
\begin{equation}
 P(n,m) = P(n-m,m) + P(n-1,m-1)
\end{equation}
\end{theorem}
\begin{proof}
\begin{equation}
\begin{split}
 & P(n-m,m) + P(n-1,m-1) 
   = [q^{n-m}] \frac{q^m}{\prod_{j=1}^m(1-q^j)} + [q^{n-1}] \frac{q^{m-1}}{\prod_{j=1}^{m-1}(1-q^j)} \\
 & = [q^n] \frac{q^{2m}}{\prod_{j=1}^m(1-q^j)} + [q^n] \frac{q^m(1-q^m)}{\prod_{j=1}^m(1-q^j)} 
  = [q^n] \frac{q^m}{\prod_{j=1}^m(1-q^j)} = P(n,m) \\
\end{split}
\end{equation}
\end{proof}
\begin{lemma}\label{mylemma}
\begin{equation}
 \sum_{k=0}^m q^k \prod_{j=k+1}^m (1-q^j) = 1
\end{equation}
\begin{proof}
The lemma is true for $m=0$, and using induction on $m$,
when it is true for $m$, then for $m+1$:
\begin{equation}
\begin{split}
 & \sum_{k=0}^{m+1} q^k \prod_{j=k+1}^{m+1} (1-q^j) = q^{m+1} + \sum_{k=0}^m q^k \prod_{j=k+1}^{m+1} (1-q^j) \\
 & = q^{m+1} + (1-q^{m+1}) \sum_{k=0}^m q^k \prod_{j=k+1}^m (1-q^j) = q^{m+1} + 1 - q^{m+1} = 1 \\
\end{split}
\end{equation}
\end{proof}
\end{lemma}
A similar lemma can be found as (5) in \cite{A83}.
Let $P(n|\textrm{\rm at most }m~\textrm{\rm parts})$ be the number of integer partitions of $n$ into at most $m$ parts \cite{AE04}.
\begin{theorem}
For integer $m\leq n$:
\begin{equation}
 P(n|\textrm{\rm at most }m~\textrm{\rm parts}) = \sum_{k=0}^m P(n,k) = P(n+m,m)
\end{equation}
\end{theorem}
\begin{proof}
Using lemma \ref{mylemma}:
\begin{equation}
\begin{split}
 & \sum_{k=0}^m P(n,k) = \sum_{k=0}^m [q^n] \frac{q^k}{\prod_{j=1}^k(1-q^j)}
  = [q^n] \sum_{k=0}^m \frac{q^k \prod_{j=k+1}^m(1-q^j)}{\prod_{j=1}^m(1-q^j)} \\
 & = [q^n] \frac{1}{\prod_{j=1}^m(1-q^j)} 
 = [q^{n+m}] \frac{q^m}{\prod_{j=1}^m(1-q^j)} 
 = P(n+m,m) \\
\end{split}
\end{equation}
\end{proof}
Taking $m=n$ it follows that $P(n)=P(2n,n)$.

\section{First Algorithm for Computing $P(n,m)$}

From the recurrence relation in theorem \ref{recur} a simple algorithm
can be derived for computing $P(n,m)$. 
Using that $P(n,m)=0$ when $n<m$,
let for each $m'$ between $1$ and $m$ an array represent $P(n',m')$
for all $n'$ between $m'$ and $n-m+m'$, starting with $P(n',1)=1$. 
Then repeated application of the recurrence relation for all $n'$ transforms the array
from $m'-1$ to $m'$, where the position in the array of $P(n',m')$
is the same as the position of $P(n'-1,m'-1)$.
When $m'=m$ is reached the array contains $P(n',m)$ for all $n'$ between $m$ and $n$,
and then the last element in the array contains $P(n,m)$.
\begin{algorithm}
  \caption{Computation of $P(n,m)$}
  \begin{algorithmic}[1]
    \Procedure{P}{$n,m$}
       \For{$p\leftarrow 0$ \textbf{to} $n-m$}
        \State $a_p\leftarrow 1$
      \EndFor
       \For{$i\leftarrow 2$ \textbf{to} $\min(m,n-m)$}
         \For{$p\leftarrow i$ \textbf{to} $n-m$}
           \State $a_p \leftarrow a_p+a_{p-i}$
         \EndFor
      \EndFor
      \State \textbf{return} $a_{n-m}$
    \EndProcedure
  \end{algorithmic}
\end{algorithm}

The number of steps $S_1(n,m)$ in this algorithm, excluding the initialization steps, is:
\begin{equation}\label{s1}
 S_1(n,m) = \frac{1}{2}(\min(m,n-m)-1)(2(n-m)-\min(m,n-m))
\end{equation}

\section{An Expression for $P(n,m)$ using $P(n)$}

An expression for $P(n,m)$ using $P(n)$ is derived,
leading to a second algorithm for computing $P(n,m)$.
Starting with theorem \ref{pnmq}:
\begin{equation}
 P(n,m) = [q^n] \frac{q^m}{\prod_{j=1}^m (1-q^j)} = [q^n] \frac{q^m \prod_{j=m+1}^{n-m}(1-q^j)}{\prod_{j=1}^n (1-q^j)} 
  = [q^n] \frac{q^m \prod_{j=1}^{n-2m}(1-q^{m+j})}{\prod_{j=1}^n (1-q^j)} 
\end{equation}
The following identity is taken from formula [3h] on page 99 in \cite{C74}:
\begin{equation}
 \prod_{j=1}^{\infty} (1-q^{m+j}) = \sum_{k=0}^{\infty} \sum_{i=0}^k (-1)^i q^{k+mi} Q(k,i)
\end{equation}
where $Q(n,m)$ is the number of partitions of $n$ into exactly $m$ distinct parts.
Because only the coefficient of $q^n$ is needed, the sum is only  needed up to $k=n-m$,
and using theorem \ref{pnq}:
\begin{equation}
\begin{split}
 P(n,m) & = [q^n] \sum_{k=0}^{n-m} \sum_{i=0}^k (-1)^i q^{k+m(i+1)} Q(k,i) \frac{1}{\prod_{j=1}^n (1-q^j)} \\
 & = \sum_{k=0}^{n-m} \sum_{i=0}^k (-1)^i Q(k,i) [q^{n-k-m(i+1)}] \frac{1}{\prod_{j=1}^n (1-q^j)} \\
 & = \sum_{k=0}^{n-m} \sum_{i=0}^k (-1)^i Q(k,i) P(n-k-m(i+1)) \\
\end{split}
\end{equation}
$Q(n,m)$ can be expressed using $P(n,m)$,
taken from page 116 in \cite{C74}, \cite{WPPQ}:
\begin{equation}
 Q(n,m) = P(n-\frac{1}{2}m(m-1),m)
\end{equation}
which results in:
\begin{equation}
 P(n,m) = \sum_{k=0}^{n-m} \sum_{i=0}^k (-1)^i P(k-\frac{1}{2}i(i-1),i) P(n-k-m(i+1))
\end{equation}
When $k=0$, only $i=0$ gives a nonzero summand because $P(0,0)=1$, leading to:
\begin{equation}
 P(n,m) = P(n-m) + \sum_{k=1}^{n-2m} \sum_{i=1}^k (-1)^i P(k-\frac{1}{2}i(i-1),i) P(n-k-m(i+1))
\end{equation}
From this formula follows:
\begin{equation}
 P(n,m) = P(n-m) \textrm{\rm~if~} m\geq \lceil n/2 \rceil
\end{equation}
Changing the order of summation gives:
\begin{equation}\label{newalgo1}
 P(n,m) = P(n-m) + \sum_{i=1}^{i_{\rm max}(n,m)} (-1)^i \sum_{k=\frac{1}{2}i(i+1)}^{n-m(i+1)} P(k-\frac{1}{2}i(i-1),i) P(n-k-m(i+1))
\end{equation}
where $i_{\rm max}(n,m)$ is given by solving $n-m(i+1)=\frac{1}{2}i(i+1)$:
\begin{equation}\label{newalgo2}
 i_{\rm max}(n,m) = \lfloor \frac{1}{2} ( \sqrt{8n+(2m-1)^2}-2m-1 ) \rfloor
\end{equation}
The two last identities are the basis for a second algorithm.

\section{Second Algorithm for Computing $P(n,m)$}

The equations (\ref{newalgo1}) and (\ref{newalgo2}) are used in the second
algorithm, where $P(0)..P(n)$ need to be computed only once,
which means they are cached in a permanent array.
In the first algorithm for each $m'$ the $P(n',m')$ were computed
for $n'$ between $m'$ and $n-m+m'$, but now they only are needed between 
$m'$ and $m'+k_{\rm max}-k_{\rm min}$.

\begin{algorithm}
  \caption{Computation of $P(n,m)$}
  \begin{algorithmic}[1]
    \Procedure{P}{$n,m$}
       \For{$p\leftarrow 0$ \textbf{to} $n-m$}
        \State $a_p\leftarrow 1$
       \EndFor
       \State $x\leftarrow P(n-m)$
       \For{$k\leftarrow 1$ \textbf{to} $n-2m$}
         \State $x \leftarrow x - P(n-2m-k)$
       \EndFor
       \For{$i\leftarrow 2$ \textbf{to} $\lfloor \frac{1}{2}(\sqrt{8n+(2m-1)^2}-2m-1)\rfloor$}
         \State $k_{\rm min}\leftarrow \frac{1}{2}i(i+1)$
         \State $k_{\rm max}\leftarrow n-m(i+1)$
         \For{$p\leftarrow i$ \textbf{to} $k_{\rm max}-k_{\rm min}$}
           \State $a_p \leftarrow a_p+a_{p-i}$
         \EndFor
         \For{$k\leftarrow k_{\rm min}$ \textbf{to} $k_{\rm max}$}
           \State $x \leftarrow x + (-1)^i a_{k-k_{\rm min}} P(k_{\rm max}-k)$
         \EndFor
      \EndFor
      \State \textbf{return} $x$
    \EndProcedure
  \end{algorithmic}
\end{algorithm}

\begin{figure}[h]
\caption{Comparison of $S_1(n,m)$ and $S_2(n,m)$ for $n=400$}
\begin{center}
\includegraphics[width=10.5cm,height=14cm,angle=-90,trim=0 1cm 0 0,keepaspectratio]{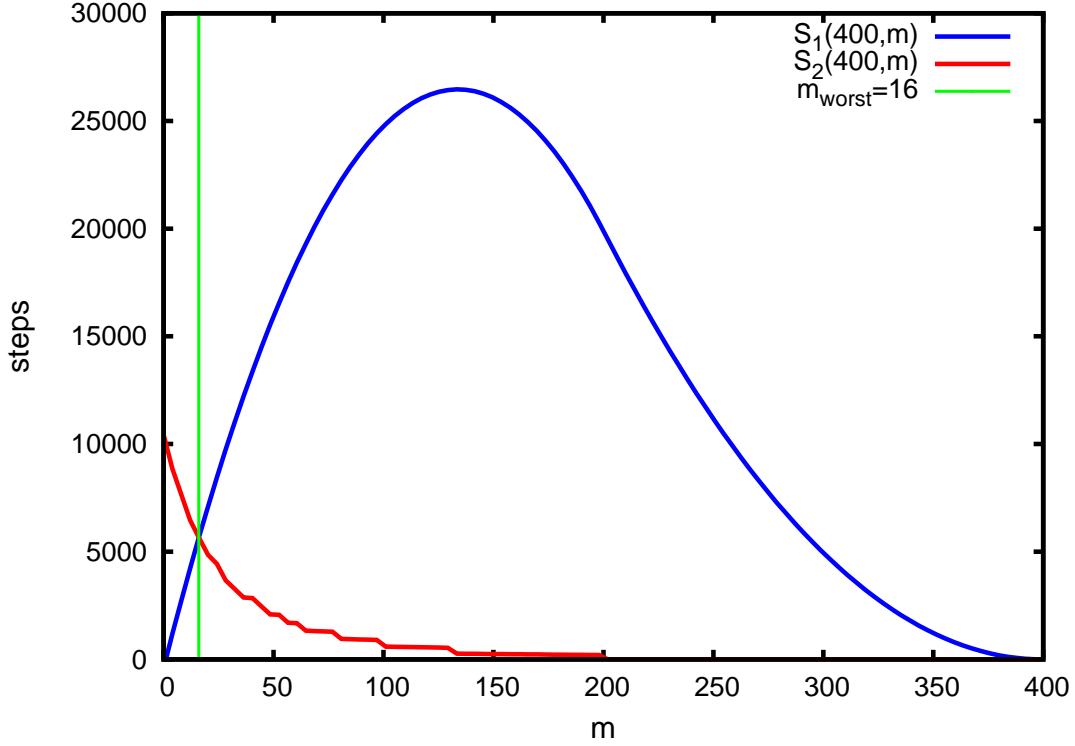}
\end{center}
\label{fig1}
\end{figure}

The number of steps $S_2(n,m)$ in this algorithm, when compared to
the first algorithm, and not including the computation of the $P(n)$
which need to be computed only once, is about:
\begin{equation}\label{s2}
 S_2(n,m) \simeq \frac{1}{2}i_{\rm max}(n,m)(2(n-m)-i_{\rm max}(n,m))
\end{equation}
In figure \ref{fig1} $S_1(n,m)$ and $S_2(n,m)$ are shown for $n=400$, from which
it is clear that for small $m$ the first algorithm and for large $m$ the second
algorithm is faster, where the cross over point is when $i_{\rm max}(n,m)=m$,
which has the solution $m=m_{\rm worst}$:
\begin{equation}
  m_{\rm worst} = \frac{1}{6}( \sqrt{24n+9} - 3 )
\end{equation}
Because in the second algorithm additional steps are needed in the loops
over $k$, in practice $m_{\rm worst}$ is a constant factor larger,
about $m_{\rm worst}\simeq 2.7\sqrt{n}$.
Substituting this $m$ in (\ref{s1}) it is clear that this combination
of these two algorithms is $O(n^{3/2})$.
Sometimes updating the cached $P(0)..P(n)$ is necessary,
but this is also $O(n^{3/2})$, see section \ref{sectlist}, 
so the total complexity remains $O(n^{3/2})$.

\section{Formulas for $P(n,m)$ when $m<=6$}

When $m\leq 6$ there are formulas for $P(n,m)$, such as $P(n,1)=1$,
$P(n,2)=\lfloor n/2\rfloor$ \cite{AE04}, and, where $[x]$ is the
nearest integer to $x$ \cite{AE04,bomfim,H85,knuth}:
\begin{equation}\label{p3}
 P(n,3) = [n^2/12]
\end{equation}
\begin{equation}
 P(n,4) = [n(2n^2+6n+9((-1)^n-1))/288]
\end{equation}
\begin{equation}
 P(n,5) = [n(n^3+10n(n+1)-15(3(-1)^n+5))/2880]
\end{equation}
\begin{equation}
 P(n,6) = [n(6n^4+135n^3+760n^2+675((-1)^n-1)n-30F(n~\textrm{mod}~6))/518400]
\end{equation}
where $F(0)=-96$, $F(1)=F(5)=629$, $F(2)=F(4)=224$ and $F(3)=309$.
For practical purposes these formulas are correct up to at least $n=10^9$.

\section{Computing a List of $P(n)$ and $Q(n)$}\label{sectlist}

For the second algorithm sometimes a list of the cached $P(0)..P(n)$ must be computed.
The first algorithm to do this is with Euler's pentagonal number
theorem \cite{A83,ewell,merca1}:
\begin{equation}
 P(n) = \delta_{n,0} -\sum_{k=1}^{\infty} (-1)^k ( P(n-k(3k-1)/2) + P(n-k(3k+1)/2) )
\end{equation}
where the sum is over all $k$ for which the argument of $P(n)$ is nonnegative.
The second algorithm with a formula of J.A. Ewell has less terms \cite{ewell,merca1}:
\begin{equation}
 P(n) = \sum_{k=0}^{\infty} P(\frac{n-k(k+1)/2}{4}) - 2 \sum_{k=1}^{\infty} (-1)^k P(n-2k^2)
\end{equation}
where the sum is over all $k$ for which the argument of $P(n)$ is a nonnegative integer.
In the computer program below this formula is about 20\% faster than the previous formula.
In the first sum the argument is an integer when $n-k(k+1)/2\equiv 0~(\textrm{mod}~4)$,
which occurs if and only if $n~\textrm{mod}~4=k(k+1)/2~\textrm{mod}~4$.
The values of $k(k+1)/2~\textrm{mod}~4$ as a function of $k$ are a repeating
pattern of $\{0,1,3,2,2,3,1,0\}$, so the $k$ begins with $k_1=\{0,1,3,2\}$ and $k_2=\{7,6,4,5\}$
depending on the value of $n~\textrm{mod}~4$, and for each term both $k_1$ and $k_2$ are increased with $8$.
Using:
\begin{equation}
 \frac{1}{2}(k+8)(k+9) - \frac{1}{2}k(k+1) = 8k+36
\end{equation}
and substituting in the right side $k=k_1+8k$ and $k=k_2+8k$ and dividing by $4$ gives
$2k_1+16k+9$ and $2k_2+16k+9$, which are the decrements of the two $P(n)$ arguments.
For efficiency the first decrement can start at $2k_1+9$ instead of $0$, 
but then for the second decrement $2k_2+9$ must be replaced by $2(k_2-k_1)$.\\
When the list $P(0)..P(n)$ is computed, a list of $Q(0)..Q(n)$ can be computed with 
another formula of J.A. Ewell \cite{ewell}:
\begin{equation}
 Q(n) = P(n) + \sum_{k=1}^{\infty} (-1)^k ( P(n-k(3k-1)) + P(n-k(3k+1))
\end{equation}
A formula of M. Merca \cite{merca2} computes a list of $Q(0)..Q(n)$ without needing $P(0)..P(n)$:
\begin{equation}
 Q(n) = s(n) -2 \sum_{k=1}^{\infty} (-1)^k Q(n-3k^2)
\end{equation}
where
\begin{equation}
 s(n) =
\begin{cases}
 1 & \text{\rm if $n=m(3m\pm 1)/2$ for some nonnegative integer $m$} \\
 0 & \text{\rm otherwise} \\
\end{cases}
\end{equation}
In the computer program below this formula is about 35\% faster than the previous formula
when all $P(0)..P(n)$ are known.
These four algorithms for computing a list of $P(n)$ and $Q(n)$ are all $O(n^{3/2})$,
and are faster than repeatedly computing isolated values of $P(n)$ and $Q(n)$
with the Hardy-Ramanujan-Rademacher formula \cite{merca1}.

\section{Computer Algebra Program}

A Mathematica\textsuperscript{\textregistered} program is given below implementing the
algorithm for $P(n,m)$ and some other related functions,
which are listed below.\\
\texttt{PartitionsPList[n]}\\
Gives a list of the $n$ numbers $P(1)..P(n)$,
where $P(n)$ is the number of partitions of $n$,
using the algorithm of J.A. Ewell \cite{ewell,merca1}.
The list is cached, and the algorithm is $O(n^{3/2})$.\\
\texttt{PartitionsQList[n]}\\
Gives a list of the $n$ numbers $Q(1)..Q(n)$,
where $Q(n)$ is the number of partitions of $n$ into distinct parts,
using the algorithm of M. Merca \cite{merca2}.
The list is cached, and the algorithm is $O(n^{3/2})$.\\
\texttt{PartitionsInPartsP[n,m]}\\
Gives the number $P(n,m)$ of partitions of $n$ into exactly $m$ parts,
with the combination of the two algorithms described in this paper,
which is $O(n^{3/2})$.\\
\texttt{PartitionsInPartsQ[n,m]}\\
Gives the number $Q(n,m)$ of partitions of $n$ into exactly $m$ distinct parts,
using the formula $Q(n,m)=P(n-m(m-1)/2,m)$, which is $O(n^{3/2})$.\\
\texttt{PartitionsInPartsPList[n]}\\
Gives a list of the $n$ numbers $P(n,1)..P(n,n)$, using the algorithm
for $P(n,m)$. This algorithm is optimized with \texttt{ListConvolve},
and is therefore $O(n^2)$.\\
\texttt{PartitionsInPartsPList[n,m]}\\
Gives a list of the $n-m+1$ numbers $P(m,m)..P(n,m)$, using the algorithm
for $P(n,m)$. This algorithm is optimized with \texttt{ListConvolve},
and is therefore $O(n^2)$.\\
\texttt{PartitionsInPartsQList[n]}\\
Gives a list of the $m_{\rm max}=\lfloor (\sqrt{8n+1}-1)/2\rfloor$ numbers
$Q(n,1)..Q(n,m_{\rm max})$, using the algorithm for $P(n,m)$.
This algorithm is $O(n^{3/2})$.\\
\texttt{PartitionsInPartsQList[n,m]}\\
Gives a list of the $n-m(m+1)/2+1$ numbers $Q(m(m+1)/2,m)..Q(n,m)$,
using the formula $Q(n,m)=P(n-m(m-1)/2,m)$.
This algorithm is $O(n^2)$.\\
\\
Below is the listing of a Mathematica\textsuperscript{\textregistered}
package which can be copied into a \texttt{PartitionsInParts.m} package file.
\newpage

\begin{verbatim}
(* ::Package:: *)

BeginPackage["PartitionsInParts`"];

PartitionsPList::usage =
"PartitionsPList[n] gives a list of the n numbers P(1)..P(n),
 the number of unrestricted partitions of 1..n.";
PartitionsQList::usage =
"PartitionsQList[n] gives a list of the n numbers Q(1)..Q(n),
 the number of partitions of 1..n into distinct parts.";
PartitionsInPartsP::usage = 
"PartitionsInPartsP[n,m] gives the number P(n,m),
 the number of partitions of n into exactly m parts.";
PartitionsInPartsQ::usage = 
"PartitionsInPartsQ[n,m] gives the number Q(n,m),
 the number of partitions of n into exactly m distinct parts.";
PartitionsInPartsPList::usage = 
"PartitionsInPartsPList[n] gives a list of the n numbers P(n,1)..P(n,n),
 the number of partitions of n into exactly 1..n parts.
PartitionsInPartsPList[n,m] gives a list of the n-m+1 numbers P(m,m)..P(n,m),
 the number of partitions of m..n into exactly m parts.";
PartitionsInPartsQList::usage = 
"PartitionsInPartsQList[n] gives a list of the numbers Q(n,1)..Q(n,mmax),
 the number of partitions of n into exactly 1..mmax distinct parts.
PartitionsInPartsQList[n,m] gives a list of the numbers Q(m(m+1)/2,m)..Q(n,m),
 the number of partitions of m(m+1)/2..n into exactly m distinct parts.";

Begin["`Private`"]; 

PartitionsPList[n_Integer?Positive]:=(partpupdateewell[n];partplist[[Range[2,n+1]]])
PartitionsQList[n_Integer?Positive]:=(partqupdatemerca[n];partqlist[[Range[2,n+1]]])
PartitionsInPartsP[n_Integer?NonNegative,m_Integer?NonNegative]:=
 If[n==m==0,1,If[m==0||n<m,0,If[n==m,1,If[m<=6,partitionsinpartsp1[n,m],
  If[m<=2.7Sqrt[n],partitionsinpartsp2[n,m,False],partitionsinpartsp3[n,m]]]]]]
PartitionsInPartsQ[n_Integer?NonNegative,m_Integer?NonNegative]:=
 If[n-m(m-1)/2<m,0,PartitionsInPartsP[n-m(m-1)/2,m]]
PartitionsInPartsPList[n_Integer?Positive]:=partitionsinpartsplist1[n]
PartitionsInPartsPList[n_Integer?NonNegative,m_Integer?NonNegative]:=
 If[m==0,PadRight[{1},n+1],If[n<m,{},
 If[m<0.21Exp[0.78Log[n]],partitionsinpartsp2[n,m,True],partitionsinpartsplist1[n,m]]]]
PartitionsInPartsQList[n_Integer?Positive]:=partitionsinpartsqlist1[n]
PartitionsInPartsQList[n_Integer?NonNegative,m_Integer?NonNegative]:=
 If[n-m(m-1)/2<m,{},PartitionsInPartsPList[n-m(m-1)/2,m]]

partplist={1};partqlist={1};
partpupdateeuler[n_]:=Block[{length=Length[partplist],result,kmax,kk},
 Which[n>=length,partplist=PadRight[partplist,n+1];
  Do[result=0;kmax=Floor[(Sqrt[24i+1]-1)/6];kk=i+1;
   Do[kk-=3k-1;result-=(-1)^k(partplist[[kk]]+partplist[[kk+k]]),{k,kmax}];
   kk-=2kmax+1;Which[kk>0,result+=(-1)^kmax partplist[[kk]]];
  partplist[[i+1]]=result,{i,length,n}]]]
partpupdateewell[n_]:=Block[{length=Length[partplist],result,ks,k1,k2,kk1,kk2},
 Which[n>=length,partplist=PadRight[partplist,n+1];
  ks={{9,11,15,13},{14,10,2,6},{0,1,6,3},{28,21,10,15}};
  Do[result=0;k1=Floor[Sqrt[i/2]];kk1=i+1;
   Do[kk1-=4k-2;result+=(-1)^k partplist[[kk1]],{k,k1}];
   result*=-2;kk2=Mod[i,4]+1;k1=ks[[1]][[kk2]];k2=ks[[2]][[kk2]];
   kk1=(i-ks[[3]][[kk2]])/4+1;kk2=(i-ks[[4]][[kk2]])/4+1;
   For[k=k1,kk2>0,kk1-=k;kk2-=k2+k;k+=16,result+=partplist[[kk1]]+partplist[[kk2]]];
   Which[kk1>0,result+=partplist[[kk1]]];
  partplist[[i+1]]=result,{i,length,n}]]]
partqupdateewell[n_]:=Block[{length=Length[partqlist],result,kmax,kk},
 Which[n>=length,partqlist=PadRight[partqlist,n+1];partpupdateewell[n];
  Do[result=partplist[[i+1]];kmax=Floor[(Sqrt[12i+1]-1)/6];kk=i+1;
   Do[kk-=6k-2;result+=(-1)^k(partplist[[kk]]+partplist[[kk+2k]]),{k,kmax}];
  kk-=4kmax+2;Which[kk>0,result-=(-1)^kmax partplist[[kk]]];
 partqlist[[i+1]]=result,{i,length,n}]]]
partqupdatemerca[n_]:=Block[{length=Length[partqlist],result,kmax,kk,tsqrt},
 Which[n>=length,partqlist=PadRight[partqlist,n+1];
  Do[result=0;kmax=Floor[Sqrt[i/3]];kk=i+1;
   Do[kk-=6k-3;result+=(-1)^k partqlist[[kk]],{k,kmax}];
   result*=-2;tsqrt=Sqrt[24i+1];
   Which[IntegerQ[(tsqrt+1)/6]||IntegerQ[(tsqrt-1)/6],result++];
  partqlist[[i+1]]=result,{i,length,n}]]]
partitionsinpartsp1[n_,m_]:=Switch[m,1,1,2,Floor[n/2],3,Round[n^2/12],
 4,Round[n(2n^2+6n+9((-1)^n-1))/288],5,Round[n(n^3+10n(n+1)-15(3(-1)^n+5))/2880],
 6,Round[n(6n^4+135n^3+760n^2+675((-1)^n-1)n-30Switch[
  Mod[n,6],0,-96,1,629,2,224,3,309,4,224,5,629])/518400]]
partitionsinpartsp2[n_,m_,all_]:=Block[{temp=ConstantArray[1,n-m+1],nmax=n-m,pmax},
 Do[pmax=(Floor[(nmax+1)/i]-1)i;
  Do[temp[[Range[p+1,p+i]]]+=temp[[Range[p-i+1,p]]],{p,i,pmax,i}];
  temp[[Range[pmax+i+1,nmax+1]]]+=temp[[Range[pmax+1,nmax-i+1]]],{i,2,Min[m,nmax]}];
 If[all,temp,temp[[nmax+1]]]]
partitionsinpartsp3[n_,m_]:=Block[{temp=ConstantArray[1,n-m+1],nmax=n-m,result,
 pmax,kmin,kmax},partpupdateewell[nmax];result=partplist[[nmax+1]];
 Do[result-=partplist[[nmax-m-k+1]],{k,nmax-m}];kmin=3;kmax=nmax-2m;
 For[i=2,kmin<=kmax,kmin+=i+1;kmax-=m;i++,nmax=kmax-kmin;pmax=(Floor[(nmax+1)/i]-1)i;
  Do[temp[[Range[p+1,p+i]]]+=temp[[Range[p-i+1,p]]],{p,i,pmax,i}];
  Which[pmax>=0,temp[[Range[pmax+i+1,nmax+1]]]+=temp[[Range[pmax+1,nmax-i+1]]]];
  Do[result+=(-1)^i temp[[k-kmin+1]]partplist[[kmax-k+1]],{k,kmin,kmax}]];result]
partitionsinpartsplist1[n_]:=Block[{temp=ConstantArray[1,n+1],
  result=ConstantArray[0,n],mmax,pmax,kmin,kmax,conv},
 mmax=Ceiling[(Sqrt[24n+9]-3)/6]+1;partpupdateewell[n-mmax];
 result[[Range[n,mmax,-1]]]=partplist[[Range[n-mmax+1]]];kmin=1;
 Do[Which[i>1,pmax=(Floor[(n+1)/i]-1)i;
  Do[temp[[Range[p+1,p+i]]]+=temp[[Range[p-i+1,p]]],{p,i,pmax,i}];
  temp[[Range[pmax+i+1,n+1]]]+=temp[[Range[pmax+1,n-i+1]]]];
  result[[i]]=temp[[n-i+1]];kmax=n-mmax(i+1)-kmin+1;
  Which[kmax>=1,conv=ListConvolve[temp[[Range[kmax]]],partplist[[Range[kmax]]],{1,1},0];
  pmax=Floor[(n-kmin)/(i+1)];
  result[[Range[pmax,mmax,-1]]]+=(-1)^i conv[[Range[kmax-(pmax-mmax)(i+1),kmax,i+1]]]];
 kmin+=i+1,{i,mmax-1}];result]
partitionsinpartsplist1[n_,m_]:=Block[{temp=ConstantArray[1,n+1],
  result=ConstantArray[0,n-m+1],mmax,pmax,kmin,kmax,conv},
 partpupdateewell[n-m];mmax=Floor[(Sqrt[8n+4m(m-1)+9]-2m-1)/2];
 result[[Range[n-m+1]]]=partplist[[Range[n-m+1]]];kmin=1;
 Do[kmax=n-m(i+1)-kmin+1;Which[i>1,pmax=(Floor[(kmax+1)/i]-1)i;
  Do[temp[[Range[p+1,p+i]]]+=temp[[Range[p-i+1,p]]],{p,i,pmax,i}];
  Which[pmax>=0,temp[[Range[pmax+i+1,kmax+1]]]+=temp[[Range[pmax+1,kmax-i+1]]]]];
  conv=ListConvolve[temp[[Range[kmax]]],partplist[[Range[kmax]]],{1,1},0];
  pmax=kmin+m(i+1);result[[Range[pmax-m+1,n-m+1]]]+=(-1)^i conv[[Range[n-pmax+1]]];
  kmin+=i+1,{i,mmax}];result]
partitionsinpartsqlist1[n_]:=Block[{temp=ConstantArray[1,n],result,mmax,nmax=n-1,pmax},
 mmax=Floor[(Sqrt[8n+1]-1)/2];result=ConstantArray[0,mmax];result[[1]]=1;
 Do[nmax-=i;pmax=(Floor[(nmax+1)/i]-1)i;
  Do[temp[[Range[p+1,p+i]]]+=temp[[Range[p-i+1,p]]],{p,i,pmax,i}];
 Which[pmax>=0,temp[[Range[pmax+i+1,nmax+1]]]+=temp[[Range[pmax+1,nmax-i+1]]]];
 result[[i]]=temp[[nmax+1]],{i,2,mmax}];result]

End[]; 
EndPackage[]; 
\end{verbatim}

\begin{figure}[h]
\caption{Timings of Mathematica\textsuperscript{\textregistered} Functions}
\begin{center}
\includegraphics[width=10.5cm,height=14cm,angle=-90,trim=0 1cm 0 0,keepaspectratio]{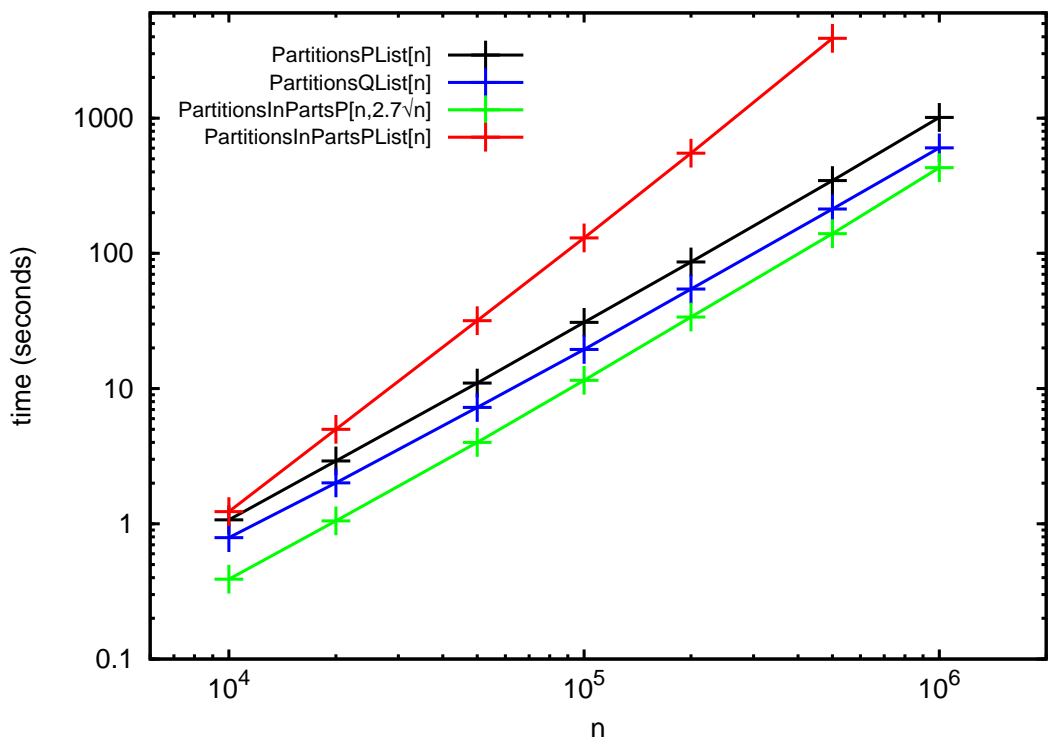}
\end{center}
\label{fig2}
\end{figure}

In figure \ref{fig2} some timings of these functions are provided
using Mathematica\textsuperscript{\textregistered} 12.3
on an Intel\textsuperscript{\textregistered} Core i7 9700K 3.60GHz with 32GB DDR4-2133 RAM.
For \texttt{PartitionsInPartsP} and \texttt{PartitionsInPartsPList} the time for updating
the list of $P(n)$, which is needed only once, is not included in the timing.
\newpage

\pdfbookmark[0]{References}{}

\end{document}